\documentclass[12pt]{amsart}
\usepackage{amsmath, amsthm, amscd, amsfonts}
\usepackage{relsize}

\newtheorem{theorem}{Theorem}[section]

\theoremstyle{definition}

\newtheorem{example}[theorem]{Example}

\theoremstyle{remark}

\newtheorem{remark}[theorem]{Remark}
\numberwithin{equation}{section}

\newfont{\kh}{msbm10}

\begin{document}
\title[Some examples of EP operators]
{Some examples of EP operators}
\author{K. Sharifi}
\address{Kamran Sharifi, \newline Faculty of Mathematical Sciences,
Shahrood University of Technology, P. O. Box 3619995161-316,
Shahrood, Iran}
\email{sharifi.kamran@gmail.com}


\begin{abstract}
We give some examples of EP and non-EP operators to show that the main results of
Mohammadzadeh Karizakia {\it et al}
[Some results about EP modular operators,
{\it Linear and Multilinear Algebra}, DOI: 10.1080/03081087.2020.1844613]
are not correct even in the case of Hilbert spaces.
\end{abstract}
\maketitle

\section{main results}
A bounded adjointable linear operator $T$ with closed range on a Hilbert C*-module (or complex
Hilbert space)  $H$ is called an $EP$ operator if $T$ and $T^*$ have
the same range. If a bounded adjointable operator $T$ does not have closed range,
then neither $Ker(T)$ nor $Ran(T)$ need to
be orthogonally complemented. For the basic theory of Hilbert C*-modules
we refer to the book \cite{LAN} and papers \cite{FR2, FR3}. Let $H$ be a Hilbert
module over an arbitrary C*-algebra of coefficients $ \mathcal{A}$. An
operator $T \in \mathcal{L}(H)$ is
called $EP$ if $Ran(T)$ and $Ran(T^*)$ have the same
closure \cite[Definition 2.1]{SHAEP}.

Mohammadzadeh Karizakia {\it et al} \cite{KARIZ1} investigate
commuting EP operators and prove that $ \overline{ Ran(T+S) } = \overline{ Ran(T) + Ran(S) }$,
when $T$ and $S$ are EP modular operators on $H$. However, the main results of this paper
is not correct even in the case of Hilbert spaces. Indeed, they have utilized some equalities and identities
that are generally not valid for operators or matrices.

\begin{example}Let bounded operators $T$ and $S$ on $ \ell_2$ be defined by
\begin{eqnarray*}
T(x_1,x_2,x_3,x_4,x_5,\ldots) &=& (x_1, x_2,x_2+x_3,x_4,x_5,\ldots) \\
S(x_1,x_2,x_3,x_4,x_5,\ldots) &=& (x_1+x_2, 0, 0, x_4,x_5,\ldots).
\end{eqnarray*}
Then
\begin{eqnarray*}
T^*(x_1,x_2,x_3,x_4,x_5,\ldots)   &=& (x_1, x_2+x_3, x_3, x_4,x_5,\ldots) \\
S^*(x_1,x_2,x_3,x_4,x_5,\ldots)   &=& (x_1, x_1, 0, x_4,x_5,\ldots)\\
(T+S)(x_1,x_2,x_3,x_4,x_5,\ldots) &=& (2x_1+x_2, x_2, x_2+x_3, x_4, x_5,\ldots)\\
(T+S)^*(x_1,x_2,x_3,x_4,x_5,\ldots) &=& (2x_1, x_1+x_2+x_3, x_3, x_4, x_5,\ldots).
\end{eqnarray*}
One can easily see that $T$ and $T+S$ are EP operators and $S$ is not an EP operator, and so the part ``$\Leftarrow$" in \cite[Theorem 2.7]{KARIZ1} is
not correct.
\end{example}
\begin{example}Let bounded operators $T$ and $S$ on $ \ell_2$ be defined by
\begin{eqnarray*}
T(x_1,x_2,x_3,x_4,x_5,\ldots) &=& (x_1-x_3, 0, x_3, x_4,x_5,\ldots) \\
S(x_1,x_2,x_3,x_4,x_5,\ldots) &=& (x_3-x_1, 0,  x_3, x_4,x_5,\ldots).
\end{eqnarray*}
Then
\begin{eqnarray*}
T^*(x_1,x_2,x_3,x_4,x_5,\ldots)   &=& (x_1, 0, x_3-x_1, x_4,x_5,\ldots) \\
S^*(x_1,x_2,x_3,x_4,x_5,\ldots)   &=& (-x_1, 0, x_3+x_1, x_4,x_5,\ldots)\\
(T+S)(x_1,x_2,x_3,x_4,x_5,\ldots) &=& (0, 0,x_3, x_4,x_5,\ldots)\\
(TT^*+SS^*)(x_1,x_2,x_3,x_4,x_5,\ldots) &=& (4x_1, 0, 2x_3, x_4,x_5,\ldots).
\end{eqnarray*}
One can easily see that $T$ and $S$ are EP operators and $S$ is not an EP operator and
\begin{eqnarray*}
\overline{ Ran(T+S) } & \neq & \overline{ Ran(TT^*) + Ran(SS^*) },\\
\overline{ Ran(T+S) } & \neq & \overline{ Ran(T) + Ran(S) },
\end{eqnarray*}
that is, the parts (iii), (iv), (v) and (vi) in \cite[Theorem 2.10]{KARIZ1} are
not correct.
\end{example}
\begin{remark}
Some mistakes of this paper based on the following gaps:
\begin{itemize}
\item  It is known that
$$ \bigoplus_{n=1} ^{ \infty} H= \{ x=(x(n)) \in \prod_n H:
~ \sum_n \langle x(n),x(n)  \rangle ~ {\rm converges~ in ~the ~ norm~of~ }  \mathcal{A} \}$$
is a Hilbert $\mathcal{A}$-module. In the proof of Theorem 2.10, the matrix operators $B'$ and $C'$ take their
values in  $\bigoplus_{n=1} ^{ \infty} H$,
and so one need to check that the matrix operators are well defined
and $\overline{Ran(B')} = \overline{Ran(C')}$.

\item Let $T$ and $S$ be bounded adjointable operators on $H$ and let
$$A= \begin{bmatrix}  T & 0 \\
S & 0 \end{bmatrix} \in \mathcal{L}(H \oplus H).
$$
The authors have applied the equality  $ \overline{ Ran(T)} \oplus \overline{ Ran(S)}  = \overline{ Ran(A)}$
several times to prove Theorem 2.7, Theorem 2.9 and Theorem 2.10. One should be aware that the equality is not valid even
for metrics.
\end{itemize}
One can consider various types of Hilbert modules and C*-algebras to reinvestigate the range equalities of
the paper \cite{KARIZ1}. In this regards, the papers \cite{FR2, FR3, FR4, SHAclosed} might be useful.
\end{remark}

\end{document}